# Regular factors of regular graphs from eigenvalues

Hongliang Lu

Center for Combinatorics, LPMC

Nankai University, Tianjin, China

**Abstract**

Let $m$ and $r$ be two integers. Let $G$ be a connected $r$-regular graph of order $n$ and $k$ an integer depending on $m$ and $r$. For even $kn$, we find a best upper bound (in terms of $r$ and $m$) on the third largest eigenvalue that is sufficient to guarantee that $G$ has a $k$-factor. When $nk$ is odd, we give a best upper bound (in terms of $r$ and $m$) on the second largest eigenvalue that is sufficient to guarantee that $G$ is $k$-critical.

## 1  Introduction

Throughout this paper, $G$ denotes a simple graph of *order* $n$ (the number of vertices) and *size* $e$ (the number of edges). For two subsets $S, T \subseteq V(G)$, let $e_G(S, T)$ denote the number of edges of $G$ joining $S$ to $T$. The eigenvalues of G are the *eigenvalues* $\lambda_i$ of its adjacent matrix $A$, indexed so that $\lambda_1 \geq \lambda_2 \cdots \geq \lambda_n$. If $G$ is $k$-regular, then it is easy to see that $\lambda_1 = k$ and also, $\lambda_2 < k$ if and only if $G$ is connected.

For a general graph $G$ and an integer $k$, a spanning subgraph $F$ such that

$$d_F(x) = k \ \ \text{for all } x \in V(G)$$

is called a *k-factor*.

Given an integer $k$ and a subgraph $H$ of $G$, we define the *deficiency* of $H$ with respect to $k$-factor as

$$def_H(G) = \sum_{v \in V} \max\{0, k - d_H(v)\}.$$

Suppose that $G$ contains no $k$-factors. Choose a spanning subgraph $F$ of $G$ such that the deficiency is minimized over all such choices. Then $F$ is called as a *k-optimal subgraph* of $G$. We call a graph $G$ *k-critical*, if $G$ contains no $k$-factors, but for any fixed vertex $x$ of $V(G)$, there exists a subgraph $H$ of $G$ such that $d_H(x) = k \pm 1$ and $d_H(y) = k$ for any vertex $y$ ($y \neq x$). Tutte [12] obtained the well-known k-Factor Theorem in 1952.

**Theorem 1.1 (Tutte [12])** *Let $k \geq 1$ be an integer and $G$ be a general graph. Then $G$ has a $k$-regular factor if and only if for all disjoint subsets $S$ and $T$ of $V(G)$,*

$$\delta(S,T) = k|S| + \sum_{x \in T} d_G(x) - k|T| - e_G(S,T) - \tau$$
$$= k|S| + \sum_{x \in T} d_{G-S}(x) - k|T| - \tau \geq 0,$$

*where $\tau$ denotes the number of components $C$, called $k$-odd components of $G - (S \cup T)$ such that $e_G(V(C), T) + k|C| \equiv 1 \pmod{2}$. Moreover, $\delta(S,T) \equiv k|V(G)| \pmod{2}$.*

Furthermore, the following well-known $k$-deficiency Theorem is due to Lovász [10] in 1970.

**Theorem 1.2 (Lovász [10])** *Let $G$ be a graph and $k$ a positive integer. Then*

$$\delta_G(S,T) = \max\{k|S| + \sum_{x \in T} d_{G-S}(x) - k|T| - q_G(S,T) \mid S, T \subseteq V(G), \text{ and } S \cap T = \emptyset\}$$

*where $q_G(S,T)$ is the number of components $C$ of $G - (S \cup T)$ such that $e(V(C), T) + k|C| \equiv 1 \pmod{2}$. Moreover, $\delta_G(S,T) \equiv k|V(G)|$.*

In [2], Brouwer and Haemers gave sufficient conditions for the existence of a 1-factor in a graph in terms of its Laplacian eigenvalues and, for a regular graph, gave an improvement in terms of the third largest adjacency eigenvalue, $\lambda_3$. Cioabă and Gregory [4] also studied relations between 1-factors and eigenvalues. Later, Cioabă, Gregory and Haemers [5] found a best upper bound on $\lambda_3$ that is sufficient to guarantee that a regular graph $G$ of order $v$ has a 1-factor when $v$ is even, and a matching of order $v - 1$ when $v$ is odd. In [11], the author studied the relation of eigenvalues and regular factors of regular graphs and obtained the following result.

**Theorem 1.3** *Let $r$ and $k$ be integers such that $1 \leq k < r$. Let $G$ be a connected $r$-regular graph with $n$ vertices and adjacency eigenvalues $r = \lambda_1 \geq \cdots \geq \lambda_n$. Let $m \geq 1$ be an integer and*

$$\lambda_3 \leq r - \frac{m-1}{r+1} + \frac{1}{(r+1)(r+2)}.$$

*Let $m^* \in \{m, m+1\}$ such that $m^* \equiv 1 \pmod{2}$. If one of the following conditions holds, then $G$ has a $k$-factor.*

(i) *$r$ is even, $k$ is odd, $|G|$ is even, and $\frac{r}{m} \leq k \leq r(1 - \frac{1}{m})$;*

(ii) $r$ is odd, $k$ is even and $k \leq r(1 - \frac{1}{m^*})$;

(iii) both $r$ and $k$ are odd and $\frac{r}{m^*} \leq k$.

In this paper, we continue to study the bound of eigenvalues of nonregular graphs and improved the bound in Theorem 1.3.

## 2 The lower bound of spectral of nonregular graphs and extremal graphs

In this section, we give the lower bound of some non-regular graphs and also gave that of the extremal graph.

**Theorem 2.1** *Let $r \geq 4$ be an integer and $m$ an even integer, where $2 \leq m \leq r+1$. Let $\mathcal{H}(r, m)$ denote the class of all connected irregular graphs with order $n \not\equiv r \pmod{2}$ and maximum degree $r$, and size $e$ with $2e \geq rn - m$. Let*

$$\rho_1(r, m) = \min_{H \in \mathcal{H}(r,m)} \lambda_1(H).$$

*Then*

$$\rho_1(r, m) = \frac{1}{2}(r - 2 + \sqrt{(r+2)^2 - 4m}). \tag{1}$$

*The extremal graph is $K_{r+1-m} + \overline{M_{m/2}}$.*

**Proof.** Let $H$ be a graph in $\mathcal{H}(r, m)$ with $\lambda_1(H) = \rho_1(r, m)$. Firstly, we prove the following claim.

*Claim 1.* $H$ has order $n$ and size $e$, where $n = r + 1$ and $2e = rn - m$.

Suppose that $2e > rn - m$. Then, since $rn - m$ is even, so $2e \geq rn - m + 2$. Because the spectra radius of a graph is at least the average degree, $\lambda_1(H) \geq \frac{2e}{n} \geq r - \frac{m-2}{r+1} > \rho_1(r, m)$. Thus $2e = rn - m$. Because $H$ has order $n$ with maximum degree $r$, we have $n \geq r + 1$. If $n > r + 1$, since $n + r$ is odd, so $n \geq r + 3$, it is straightforward to check that

$$\lambda_1(H) > \frac{2e}{n} \geq r - \frac{m}{r+3} > \rho_1(r, m),$$

a contradiction. We complete the claim.

Then by Claim 1, $H$ has order $n = r + 1$ and at least $r + 1 - m$ vertices of degree $r$. Let $G_1$ be the subgraph of $H$ induced by $n_1 = n + 1 - m$ vertices of all the vertices of degree $r$ and $G_2$ be the subgraph induced by the remaining $n_2 = m$ vertices. Also, let $G_{12}$ be the bipartite subgraph induced by the partition and let $e_{12}$ be the size of $G_{12}$. A theorem of Hamers [7] shows that eigenvalues of the quotient matrix of the partition interlace the eigenvalues of the adjacency matrix of $G$. Because each vertex in $G_1$ is adjacent to all other vertices in $H$. Here the quotient matrix is

$$Q = \begin{pmatrix} \frac{2e_1}{n_1} & \frac{e_{12}}{n_1} \\ \frac{e_{12}}{n_2} & \frac{2e_2}{n_2} \end{pmatrix} = \begin{pmatrix} r-m & m \\ r+1-m & m-2 \end{pmatrix}.$$

Applying this result to the greatest eigenvalue of $G$, we get

$$\lambda_1(H) \geq \lambda_1(Q) = \frac{1}{2}(r-2+\sqrt{(r+2)^2-4m}),$$

with the equality if and only if the partition is equitable [ [9], p. 195]; equivalently, if and only if $G_1$ and $G_2$ are regular, and $G_{12}$ is semiregular. So $G_2 = \overline{K_{m/2}}$, $G_1 = K_{r+1-m}$ and $G_{12} = K_{r+1-m,m}$. Hence $H = K_{r+1-m} + \overline{M_{m/2}}$. $\square$

**Theorem 2.2** *Let $r$ and $m$ be two integers such that $m \equiv r \pmod{2}$ and $1 \leq m \leq r+1$. Let $\mathcal{H}(r,m)$ denote the class of all connected irregular graphs with order $n \equiv r \pmod{2}$, maximum degree $r$, and size $e$ with $2e \geq rn - m$. Let*

$$\rho(r,m) = \min_{H \in \mathcal{H}(r,m)} \lambda_1(H).$$

(i) *If $m \geq 3$, then*

$$\rho(r,m) = \frac{1}{2}(r-3+\sqrt{(r+3)^2-4m}), \qquad (2)$$

*and then the extremal graph is $\overline{M_{(r+2-m)/2}} + \overline{C}$, where $C$ with order $m$ consists of disjoint cycles;*

(ii) *if $m = 1$, then $\rho(r,m)$ is the greatest root of $P(x)$, where $P(x) = x^3 - (r-2)x^2 - 2rx + (r-1)$;*

(iii) *if $m = 2$, then $\rho(r,m)$ is the greatest root of $f_1(x)$, where $f_1(x) = x^3 - (r-2)x^2 - (2r-1)x + r$.*

**Proof.** Let $H$ be a graph in $\mathcal{H}(r,m)$ with $\lambda_1(H) = \rho(r,m)$. With similar proof of Claim 1 in Theorem 2.1, we obtain the following claim.

*Claim 1.* $H$ has order $n$ and size $e$, where $n = r+2$ and $2e = rn - m$.

By Claim 1, $H$ has order $n = r+2$ and at least $r+2-m$ vertices of degree $r$. Let $G_1$ be the subgraph of $H$ induced by the $n_1 = n+2-m$ vertices of degree $r$ and $G_2$ be the subgraph induced by the remaining $n_2 = m$ vertices. Also, let $G_{12}$ be the bipartite subgraph induced by the partition and let $e_{12}$ be the size of $G_{12}$. Here the quotient matrix is

$$Q = \begin{pmatrix} \frac{2e_1}{n_1} & \frac{e_{12}}{n_1} \\ \frac{e_{12}}{n_2} & \frac{2e_2}{n_2} \end{pmatrix}.$$

Suppose that $e_{12} = t$. Then $2e_1 = (r+2-m)r - t$ and $2e_2 = rm - m - t$. Applying this result to greatest eigenvalue

$$\lambda_1(G) \geq \lambda_1(Q) = \frac{2e_1}{n_1} + \frac{2e_2}{n_2} + \sqrt{(\frac{2e_1}{n_1} - \frac{2e_2}{n_2})^2 + \frac{e_{12}^2}{n_1 n_2}}$$

$$= \frac{2r-1}{2} - \frac{(r+2)t}{2m(r+2-m)} + \sqrt{(\frac{1}{2} + \frac{t(r+2-2m)}{2m(r+2-m)})^2 + \frac{t^2}{m(r+2-m)}}.$$

Let $s = \frac{t}{m(r+2-m)}$, where $0 < s \leq 1$, then we have

$$2\lambda_1(Q) = f(s) = (2r-1) - s(r+2) + \sqrt{1 + 2s(r+2-2m) + s^2(r+2)^2}.$$

For $s > 0$, since

$$f'(s) = -(r+2) + \frac{(r+2-2m) + s(r+2)^2}{\sqrt{1 + 2s(r+2-2m) + s^2(r+2)^2}} < 0.$$

Then $0 < t \leq m(r+1-m)$, so we have

$$2\lambda_1(Q) \geq f(1) = (r-3) + \sqrt{1 + 2(r+2-2m) + (r+2)^2}$$
$$= (r-3) + \sqrt{(r+3)^2 - 4m}.$$

Hence

$$\lambda_1(H) \geq \lambda_1(Q) \geq \frac{1}{2}(r-3) + \frac{1}{2}\sqrt{(r+3)^2 - 4m},$$

with equality if and only if $t = m(r+2-m)$, both $G_1$ and $G_2$ are regular and $G_{12}$ is semiregular. So $\overline{G_1}$ is a perfect matching and $\overline{G_2}$ is a 2-regular graph.

*Case 2.* $m = 1$.

Then $r$ is odd and $n = r + 2$. So $H$ contains one vertex of degree $r - 1$, say $v$ and the rest vertices have degree $r$. Hence $\overline{H} = K_{1,2} \cup M_{(r-1)/2}$. Partition the vertex of $V(\overline{H})$ into three parts: the two endpoints of $K_{1,2}$; the internal vertex of $K_{1,2}$; the $(r-1)$ vertices of $M_{(r-1)/2}$. This is an equitable partition of $H$ with quotient matrix

$$Q = \begin{pmatrix} 0 & 0 & r-1 \\ 0 & 1 & r-1 \\ 1 & 2 & r-3 \end{pmatrix}.$$

The characteristic polynomial of the quotient matrix is

$$P(x) = x^3 - (r-2)x^2 - 2rx + (r-1).$$

Since the partition is equitable, so $\lambda_1(H) = \lambda_1(Q)$ and $\lambda_1(H)$ is a root of $P(x)$.

*Case 3.* $m = 2$.

Then $r$ is even. Let $G \in \mathcal{H}(r,m)$ be the graph with order $r+2$ and size $e = (r(r+2)-2)/2$. We discuss three subcases.

*Subcase 3.1.* $G$ has two nonadjacent vertices of degree $r-1$.

Then $\overline{G} = P_4 \cup M_{(r-2)/2}$. Partition the vertex of $V(\overline{G})$ into three parts: the two endpoints of $P_4$; the two internal vertices of $P_4$; the $(r-2)$ vertices of $M_{(r-1)/2}$. This is an equitable partition of $G$ with quotient matrix

$$Q_1 = \begin{pmatrix} 1 & 1 & r-2 \\ 0 & 1 & r-2 \\ 2 & 2 & r-4 \end{pmatrix}.$$

The characteristic polynomial of the quotient matrix is

$$f_1(x) = x^3 - (r-2)x^2 - (2r-1)x + r.$$

*Subcase 3.2.* $G$ has two adjacent vertices of degree $r-1$.

Then $\overline{G} = 2P_3's \cup M_{(r-4)/2}$. Still partition the vertex of $V(\overline{G})$ into three parts: the four endpoints of $P_3$; the two internal vertices of two $P_3$; the $(r-4)$ vertices of $M_{(r-1)/2}$. This is an equitable partition of $G$ with quotient matrix

$$Q_3 = \begin{pmatrix} 3 & 1 & r-4 \\ 2 & 1 & r-4 \\ 4 & 2 & r-6 \end{pmatrix}.$$

The characteristic polynomial of the quotient matrix is

$$f_2(x) = x^3 - (r-2)x^2 - (2r-1)x + r - 2.$$

*Subcase 3.3.* $G$ has one vertex of degree $r-2$.

Then $\overline{G} = K_{1,3} \cup M_{(r-2)/2}$. Partition the vertex of $V(\overline{G})$ into three parts: the center vertex of $K_{1,3}$; the three endpoints of $K_{1,3}$; the $(r-2)$ vertices of $M_{(r-2)/2}$. This is an equitable partition of $G$ with quotient matrix

$$Q_2 = \begin{pmatrix} 0 & 0 & r-2 \\ 0 & 2 & r-2 \\ 1 & 3 & r-4 \end{pmatrix}.$$

The characteristic polynomial of the quotient matrix is

$$f_3(x) = x^3 - (r-2)x^2 - 2rx + 2(r-2).$$

Note that $\lambda_1(Q_1) < \lambda_1(Q_2) < \lambda_1(Q_3)$. We have $\rho_2(r,m) = \lambda_1(Q_1)$. So $\overline{H} = P_4 \cup M_{(r-2)/2}$. Hence $\lambda_1(H)$ is a root of $f_1(x) = 0$. This completes the proof. $\square$

## 3 Regular factor of regular graph

**Lemma 3.1** *Let $r$ and $k$ be integers such that $1 \leq k < r$. Let $G$ be a connected $r$-regular graph with $n$ vertices. Let $m$ be an integer such that $m^* \in \{m, m+1\}$ and $m^* \equiv 1 \pmod{2}$. Suppose that one of the following conditions holds*

*(i) $r$ is even, $k$ is odd, $|G|$ is even, and $\frac{r}{m} \leq k \leq r(1 - \frac{1}{m})$;*

*(ii) $r$ is odd, $k$ is even and $k \leq r(1 - \frac{1}{m^*})$;*

*(iii) both $r$ and $k$ are odd and $\frac{r}{m^*} \leq k$.*

*If $G$ contains no a $k$-factor and is not $k$-critical, then $G$ contains $def(G)+1$ vertex disjoint induced subgraph $H_1, H_2, \ldots, H_{def(G)+1}$ such that $2e(H_i) \geq r|V(H_i)| - (m-1)$ for $i = 1, 2, \ldots, def(G) + 1$.*

**Proof.** Suppose that the result doesn't hold. Let $\theta = k/r$. Since $G$ is not $k$-critical and contains no $k$-factor, so by Theorem 1.2, there exist two disjoint subsets $S$ and $T$ of $V(G)$ such that $S \cup T \neq \emptyset$ and $\delta(S,T) = -def(G) \leq -1$. Let $C_1, \ldots, C_\tau$ be the $k$-odd components of $G - (S \cup T)$. We have

$$-def(G) = \delta(S,T) = k|S| + \sum_{x \in T} d_G(x) - k|T| - e_G(S,T) - \tau. \qquad (3)$$

Claim 1. $\tau \geq def(G) + 1$.

Otherwise, let $\tau \leq def(G)$. Then we have

$$0 \geq k|S| + \sum_{x \in T} d_{G-S}(x) - k|T|. \qquad (4)$$

So we have $|S| \leq |T|$, and equality holds only if $\sum_{x \in T} d_{G-S}(x) = 0$. Since $G$ is $r$-regular, so we have

$$r|S| \geq e_G(S,T) = r|T| - \sum_{x \in T} d_{G-S}(x). \tag{5}$$

By (4) and (5), we have

$$(r-k)(|T| - |S|) \leq 0.$$

Hence $|T| = |S|$ and $\sum_{x \in T} d_{G-S}(x) = 0$. Since $G$ is connected, so we have $\tau = def(G) > 0$. Since $G$ is connected, then $e_G(C_i, S \cup T) > 0$ and so $e_G(C_1, S) > 0$. Note that $G$ is $r$-regular, then we have $r|S| \geq r|T| - \sum_{x \in T} d_{G-S}(x) + e(C_i, S)$, a contradiction. We complete the claim.

By the hypothesis, without loss of generality, we can say $e(S \cup T, C_i) \geq m$ for $i = 1, \ldots, \tau - def(G)$. Then $0 < \theta < 1$, and we have

$$-def(G) \geq \delta(S,T) = k|S| + \sum_{x \in T} d_G(x) - k|T| - e_G(S,T) - \tau$$

$$= k|S| + (r-k)|T| - e_G(S,T) - \tau$$

$$= \theta r|S| + (1-\theta)r|T| - e_G(S,T) - \tau$$

$$= \theta \sum_{x \in S} d_G(x) + (1-\theta) \sum_{x \in T} d_G(x) - e_G(S,T) - \tau$$

$$\geq \theta(e_G(S,T) + \sum_{i=1}^{m} e_G(S, C_i)) + (1-\theta)(e_G(S,T) + \sum_{i=1}^{m} e_G(T, C_i)) - e_G(S,T) - \tau$$

$$= \sum_{i=1}^{m} (\theta e_G(S, C_i) + (1-\theta)e_G(T, C_i) - 1).$$

Since $G$ is connected, so we have $\theta e_G(S, C_i) + (1-\theta)e_G(T, C_i) > 0$ for $1 \leq i \leq \tau$. Hence it suffices to show that for every $C = C_i$, $1 \leq i \leq \tau - def(G)$,

$$\theta e_G(S, C_i) + (1-\theta)e_G(T, C_i) \geq 1. \tag{6}$$

Since $C$ is a $k$-odd component of $G - (S \cup T)$, we have

$$k|C| + e_G(T, C) \equiv 1 \pmod{2}. \tag{7}$$

Moreover, since $r|C| = e_G(S \cup T, C) + 2|E(C)|$, then we have

$$r|C| \equiv e_G(S \cup T, C) \pmod{2}. \tag{8}$$

It is obvious that the two inequalities $e_G(S, C) \geq 1$ and $e_G(T, C) \geq 1$ implies

$$\theta e_G(S, C) + (1-\theta)e_G(T, C) \geq \theta + (1-\theta) = 1.$$

Hence we may assume $e_G(S,C) = 0$ or $e_G(T,C) = 0$. We consider two cases.

First we consider (i). If $e_G(S,C) = 0$, since $1 \le k \le r(1 - \frac{1}{m})$, then $\theta \le 1 - \frac{1}{m}$ and so $1 \le (1-\theta)m$. Note that $e(T,C) \ge m$, so we have

$$(1-\theta)e_G(T,C) \ge (1-\theta)m \ge 1.$$

If $e_G(T,C) = 0$, since $k \ge r/m$, so $m\theta \ge 1$. Hence we obtain

$$\theta e_G(S,C) \ge m\theta \ge 1.$$

In order to prove that (ii) implies the claim, it suffices to show that (6) holds under the assumption that $e_G(S,C)$ or $e_G(T,C) = 0$. If $e_G(S,C) = 0$, then by (7), we have $e_G(T,C) \equiv 1 \pmod 2$. Hence $e_G(T,C) \ge m^*$, and thus

$$(1-\theta)e_G(T,C) \ge (1-\theta)m^* \ge 1.$$

If $e_G(T,C) = 0$, then by (8), we have $k|C| \equiv 1 \pmod 2$, which contradicts the assumption that $k$ is even.

We next consider (iii), i.e., we assume that both $r$ and $k$ are odd and $\frac{r}{m^*} \le k$. If $e_G(S,C) = 0$, then by (7) and (8), we have

$$|C| + e_G(T,C) \equiv 1 \pmod 2 \text{ and } |C| \equiv e_G(T,C) \pmod 2.$$

This is a contradiction. If $e_G(T,C) = 0$, then by (7) and (8), we have

$$|C| \equiv 1 \pmod 2 \text{ and } |C| \equiv e_G(S,C) \pmod 2,$$

which implies $e_G(S,C) \ge m^*$. Thus

$$\theta e_G(S,C) \ge \theta m^* \ge 1.$$

So we have

$$-def(G) \ge \delta(S,T) > -def(G),$$

a contradiction. This completes the proof. □

**Theorem 3.2** *Suppose that $r$ is even, $k$ is odd. Let $G$ be a connected $r$-regular graph with order $n$. Let $m \ge 3$ be an integer and $m_0 \in \{m, m-1\}$ be an odd integer. Suppose that $\frac{r}{m} \le k \le r(1 - \frac{1}{m})$.*

(i) *If $n$ is odd and $\lambda_2(G) < \rho_1(r, m_0 - 1)$, then $G$ is $k$-critical;*

(ii) if $n$ is even and $\lambda_3(G) < \rho_1(r, m_0 - 1)$, then $G$ contains $k$-factor.

**Proof.** Firstly, we prove (i). Conversely, suppose that $G$ is not $k$-critical. By Lemma 3.1, $G$ contains two vertex disjoint induced subgraphs $H_1$ and $H_2$ such that $2e(H_i) \geq rn_i - (m-1)$, where $n_i = |V(H_i)|$ for $i = 1, 2$. Hence we have $2e(H_i) \geq rn_i - (m_0 - 1)$. So by Interlacing Theorem, we have

$$\lambda_2(G) \geq \min\{\lambda_1(H_1), \lambda_1(H_2)\}$$
$$\geq \min\{\rho_1(r, m_0 - 1), \rho_2(r, m_0 - 1)\} = \rho_1(r, m_0 - 1).$$

So we have $\lambda_2(G) \geq \rho_1(r, m_0 - 1)$, a contradiction.

Now we prove (ii). Suppose that $G$ contains no a $k$-factor. Then we have $def(G) \geq 2$. So by Lemma 3.1, $G$ contains three vertex disjoint induced subgraphs $H_1$, $H_2$ and $H_3$ such that $2e(H_i) \geq rn_i - (m-1)$, where $n_i = |V(H_i)|$ for $i = 1, 2, 3$. Since $r$ is even, so $2e(H_i) \geq rn_i - (m_0 - 1)$ for $i = 1, 2, 3$. So by Interlacing Theorem, we have

$$\lambda_3(G) \geq \min\{\lambda_1(H_1), \lambda_1(H_2), \lambda_1(H_3)\}$$
$$\geq \min\{\rho_1(r, m_0 - 1), \rho_2(r, m_0 - 1)\} = \rho_1(r, m_0 - 1),$$

a contradiction. We complete the proof. □

**Theorem 3.3** *Let $r$ and $k$ be two integers. Let $m$ be an integer such that $m^* \in \{m, m+1\}$ and $m^* \equiv 1 \pmod 2$. Let $G$ be a connected $r$-regular graph with order $n$. Suppose that*

$$\lambda_3(G) < \begin{cases} \rho_1(r, m-1) & \text{if } m \text{ is odd,} \\ \rho_2(r, m-1) & \text{if } m \text{ is even.} \end{cases}$$

*If one of the following conditions holds, then $G$ has a $k$-factor.*

*(i) $r$ is odd, $k$ is even and $k \leq r(1 - \frac{1}{m^*})$;*

*(ii) both $r$ and $k$ are odd and $\frac{r}{m^*} \leq k$.*

**Proof.** Suppose that $G$ contains no a $k$-factor. By Lemma 3.1, $G$ contains three vertex disjoint induced subgraph $H_1, H_2, H_3$ such that $2e(H_i) \geq r|V(H_i)| - (m-1)$ for $i = 1, 2, 3$. Firstly, let $m$ be odd. By Interlacing Theorem we have

$$\lambda_3(G) \geq \min_{1 \leq i \leq 3} \lambda_1(H_i) \geq \min\{\rho_1(r, m-1), \rho_2(r, m-2)\} = \rho_1(r, m-1).$$

So we have $\lambda_3(G) \geq \rho_2(r, m-1)$, a contradiction.

Secondly, let $m$ be even. By Interlacing Theorem we have

$$\lambda_3(G) \geq \min_{1 \leq i \leq 3} \lambda_1(H_i) \geq \min\{\rho_1(r, m-2), \rho_2(r, m-1)\} = \rho_2(r, m-1),$$

a contradiction. We complete the proof. □